\documentclass[11pt,a4paper]{article}
\usepackage{amsbsy,amsfonts,amssymb,amsmath,amsthm}
\usepackage{tikz}

\theoremstyle{remark}

\title{Counting distinct dimer hex tilings}
\author{Peter Taylor}
\date{Preprint 2015-06-17\footnote{A minor update was made on 2016-02-17 to incorporate some corrections from N. J. A. Sloane.}}

\begin{document}
\maketitle

The combinatorics of tilings of a hexagon of integer side-length $n$ by $120^\circ-60^\circ$ diamonds of side-length $1$ has a long history, both directly (as a problem of interest in thermodynamic models) and indirectly (through the equivalence to plane partitions \cite{David89}). Formulae as products of factorials have been conjectured and, one by one, proven for the number of such tilings under each of the symmetries of the hexagon. However, when this note was written the entry for the number of \textit{distinct} such tilings in the Online Encyclopedia of Integer Sequences (OEIS) consisted of little more than a table for $0 \le n \le 4$ and a brief discussion of those values \cite{A066931}. The aim of this note is to pull together the relevant facts and allow the entry to be improved.

The symmetry group of the hexagon is the dihedral group of order $12$, $D_{12}$, which has the presentation $\langle f, r \mid f^2 = r^6 = (fr)^2 = e\rangle$ where, for our geometric purposes, $f$ (for ``flip'') represents a reflection around a diameter and $r$ (for ``rotation'') represents a rotation by $60^\circ$. At a group-theoretic level it doesn't matter whether we take the reflection to be around a diameter that passes through two vertices of the hexagon, or through the midpoints of two edges; but for the purposes of this note it is taken to pass through two vertices. Unfortunately, $fr$ then corresponds to Stanley's complementation tranform, so a rewrite may be in order to unify the two presentations.

The subgroup poset of $D_{12}$ contains $16$ subgroups, and the poset can be expressed entirely in terms of a handful of parameterised maximal chains:

\begin{tabular}{ll}
$\langle e\rangle < \langle r^a \rangle < \langle r \rangle < \langle f, r \rangle$ & $a \in \{2, 3\}$ \\
$\langle e\rangle < \langle r^a \rangle < \langle fr^b, r^a \rangle < \langle f, r \rangle$ & $a \in \{2, 3\}$, $b \in \{0, 1\}$ \\
$\langle e\rangle < \langle fr^a \rangle < \langle fr^{a \mod b}, r^b \rangle < \langle f, r \rangle$ & $a \in \{0, 1\}$, $b \in \{2, 3\}$
\end{tabular}

If we denote the number of tilings which have at least the symmetries of $g \le \langle f, r\rangle$ by $\#c(g)$ then the number of tilings having exactly the symmetries of $g$ is given by the M\"obius function of the incidence algebra of this poset: $\sum_{h \le \langle f, r\rangle} \mu(g, h) \#c(h)$. We then normalise each by its index in $\langle f, r\rangle$, since we want the number of tilings corresponding to a single coset. Thus the desired value is $$\textrm{A066931}(n) = \sum_{g \le \langle f, r\rangle} \frac{1}{[D_{12} : g]} \sum_{h \le \langle f, r\rangle} \mu(g, h) \#c(h)$$

At this point, it is worth noting that although there are $16$ subgroups, there are only $10$ conjugacy groups because the $9$ subgroups which contain only reflections each have three choices of axis / axes. Thus there is no disagreement with Stanley's list of 10 symmetries of plane partitions \cite{Stanley86}.

Once we take that into account and work through the numbers, we find that some of the weights cancel out, leaving $$\textrm{A066931}(n) = \frac{\#c\langle e\rangle + 3\#c\langle f\rangle + 3\#c\langle fr\rangle +\#c\langle r^3\rangle + 2\#c\langle r^2\rangle + 2\#c\langle r\rangle}{12}$$

$\#c\langle e\rangle$ is the most studied of the cases, being simply the number of tilings (Stanley's case 1). It is A008793 in OEIS, and has formula $\#c\langle e\rangle(n) = \prod_{i=0}^{n-1} \frac{i!(i+2n)!}{(i+n)!^2}$ \cite{MacMahon16}.

$\#c\langle f\rangle$ is A049505, the ``number of symmetric plane partitions in n-cube'' (Stanley's case 2), with formula $\#c\langle f\rangle(n) = \prod_{i=0}^{n-1} \frac{(2i)!(i+2n)!}{(2i+n)!(i+n)!}$ (with a change of index from the OEIS entry for consistency with $\#c\langle e\rangle$) \cite{Andrews78}.

$\#c\langle fr\rangle$ is A181119, ``transpose-complementary plane partitions'' (Stanley's case 6), but there's a slight catch. It is non-zero only for even $n$; there is a simple geometric explanation. The axis of symmetry bisects $2n$ triangles, which must pair up along the axis because otherwise they would be forced to pair with two triangles to maintain the symmetry. But the pairs start at the centre of the hexagon, so if $n$ is odd there are two triangles on opposite edges which can't be paired. (See figure 1). The upshot is that the given formula must be inflated as $\#c\langle fr\rangle(2m) = \binom{3m-1}{m} \prod_{i=0}^{2m-3} \frac{(i+4m)!(2i+2)!}{(2i+2m+2)!(i+2m)!}$ \cite{Proctor88}.

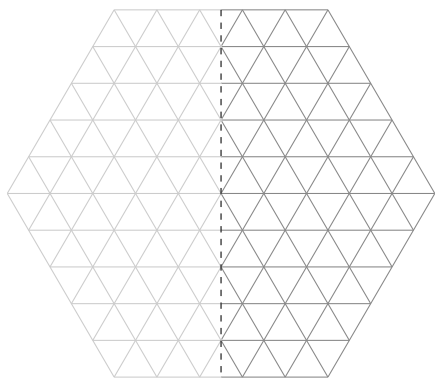
\begin{figure}[htbp]
	\centering
	\definecolor{c888}{RGB}{136,136,136}
	\definecolor{cccc}{RGB}{204,204,204}
	\begin{tikzpicture}[y=0.80pt, x=0.80pt, yscale=-1.000000, xscale=1.000000, inner sep=0pt, outer sep=0pt]
	\begin{scope}[shift={(100.0,100.0)}]
	\begin{scope}[draw=c888]
		\clip(0,-100) rectangle (100,100);
		\path[draw]
			(-50.0000,-86.6025) -- (50.0000,-86.6025) (60.0000,-69.2820) -- (-60.0000,-69.2820)
			(-70.0000,-51.9615) -- (70.0000,-51.9615) (80.0000,-34.6410) -- (-80.0000,-34.6410)
			(-90.0000,-17.3205) -- (90.0000,-17.3205) (100.0000,0.0000) -- (-100.0000,0.0000)
			(-90.0000,17.3205) -- (90.0000,17.3205) (80.0000,34.6410) -- (-80.0000,34.6410)
			(-70.0000,51.9615) -- (70.0000,51.9615) (60.0000,69.2820) -- (-60.0000,69.2820)
			(-50.0000,86.6026) -- (50.0000,86.6026);
		\begin{scope}[rotate=120.0]
			\path[draw]
				(-50.0000,-86.6025) -- (50.0000,-86.6025) (60.0000,-69.2820) -- (-60.0000,-69.2820)
				(-70.0000,-51.9615) -- (70.0000,-51.9615) (80.0000,-34.6410) -- (-80.0000,-34.6410)
				(-90.0000,-17.3205) -- (90.0000,-17.3205) (100.0000,0.0000) -- (-100.0000,0.0000)
				(-90.0000,17.3205) -- (90.0000,17.3205) (80.0000,34.6410) -- (-80.0000,34.6410)
				(-70.0000,51.9615) -- (70.0000,51.9615) (60.0000,69.2820) -- (-60.0000,69.2820)
				(-50.0000,86.6026) -- (50.0000,86.6026);
		\end{scope}
		\begin{scope}[rotate=240.0]
			\path[draw]
				(-50.0000,-86.6025) -- (50.0000,-86.6025) (60.0000,-69.2820) -- (-60.0000,-69.2820)
				(-70.0000,-51.9615) -- (70.0000,-51.9615) (80.0000,-34.6410) -- (-80.0000,-34.6410)
				(-90.0000,-17.3205) -- (90.0000,-17.3205) (100.0000,0.0000) -- (-100.0000,0.0000)
				(-90.0000,17.3205) -- (90.0000,17.3205) (80.0000,34.6410) -- (-80.0000,34.6410)
				(-70.0000,51.9615) -- (70.0000,51.9615) (60.0000,69.2820) -- (-60.0000,69.2820)
				(-50.0000,86.6026) -- (50.0000,86.6026);
		\end{scope}
	\end{scope}
	\begin{scope}[draw=cccc]
		\clip(-100,-100) rectangle (0,100);
		\path[draw]
				(-50.0000,-86.6025) -- (50.0000,-86.6025) (60.0000,-69.2820) -- (-60.0000,-69.2820)
				(-70.0000,-51.9615) -- (70.0000,-51.9615) (80.0000,-34.6410) -- (-80.0000,-34.6410)
				(-90.0000,-17.3205) -- (90.0000,-17.3205) (100.0000,0.0000) -- (-100.0000,0.0000)
				(-90.0000,17.3205) -- (90.0000,17.3205) (80.0000,34.6410) -- (-80.0000,34.6410)
				(-70.0000,51.9615) -- (70.0000,51.9615) (60.0000,69.2820) -- (-60.0000,69.2820)
				(-50.0000,86.6026) -- (50.0000,86.6026);
		\begin{scope}[rotate=120.0]
			\path[draw]
				(-50.0000,-86.6025) -- (50.0000,-86.6025) (60.0000,-69.2820) -- (-60.0000,-69.2820)
				(-70.0000,-51.9615) -- (70.0000,-51.9615) (80.0000,-34.6410) -- (-80.0000,-34.6410)
				(-90.0000,-17.3205) -- (90.0000,-17.3205) (100.0000,0.0000) -- (-100.0000,0.0000)
				(-90.0000,17.3205) -- (90.0000,17.3205) (80.0000,34.6410) -- (-80.0000,34.6410)
				(-70.0000,51.9615) -- (70.0000,51.9615) (60.0000,69.2820) -- (-60.0000,69.2820)
				(-50.0000,86.6026) -- (50.0000,86.6026);
		\end{scope}
		\begin{scope}[rotate=240.0]
			\path[draw]
				(-50.0000,-86.6025) -- (50.0000,-86.6025) (60.0000,-69.2820) -- (-60.0000,-69.2820)
				(-70.0000,-51.9615) -- (70.0000,-51.9615) (80.0000,-34.6410) -- (-80.0000,-34.6410)
				(-90.0000,-17.3205) -- (90.0000,-17.3205) (100.0000,0.0000) -- (-100.0000,0.0000)
				(-90.0000,17.3205) -- (90.0000,17.3205) (80.0000,34.6410) -- (-80.0000,34.6410)
				(-70.0000,51.9615) -- (70.0000,51.9615) (60.0000,69.2820) -- (-60.0000,69.2820)
				(-50.0000,86.6026) -- (50.0000,86.6026);
		\end{scope}
	\end{scope}
	\path[draw=black,dash pattern=on 2.40pt off 2.40pt] (0.0000,-86.6025) -- (0.0000,86.6025);
	\end{scope}

	\end{tikzpicture}

	\caption{for $\langle fr\rangle$ symmetry, the small triangles which the axis bisects must pair up along the axis.}
\end{figure}

$\#c\langle r^3\rangle$ is, astonishingly, not in OEIS!\footnote{Update: it is now A259049.} It corresponds to Stanley's case 5 (self-complementary plane partitions), and is non-zero only for even $n$, for a slightly different reason. If we divide the hexagon into six large triangles by drawing the diameters through its vertices, each of these large triangles has $\frac{n(n+1)}{2}$ small triangles in one orientation and $\frac{(n-1)n}{2}$ small triangles in another orientation. For $\langle r^3\rangle$ we can consider a half-hexagon formed from three of the large triangles: the symmetry will give us the tiling of the other half. This half-hexagon has a surplus of $n$ triangles in one orientation; therefore there must be $n$ triangles along the diameter which pair across into the other half-hexagon. (See figure 2). By the rotational symmetry, these $n$ triangles must be paired up, but if $n$ is odd this is impossible. For even $n$ we have $\#c\langle r^3\rangle(2m) = \#c\langle e\rangle(m)^2$ \cite{Stanley86}.

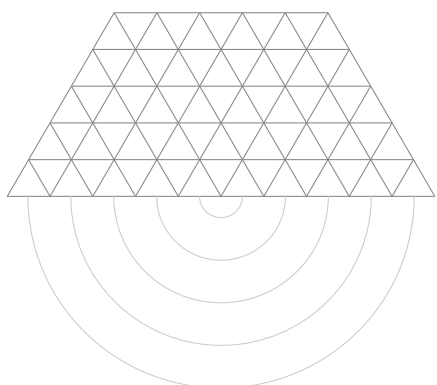
\begin{figure}[htbp]
	\centering
	\definecolor{c888}{RGB}{136,136,136}
	\definecolor{cccc}{RGB}{204,204,204}
	\begin{tikzpicture}[y=0.80pt, x=0.80pt, yscale=-1.000000, xscale=1.000000, inner sep=0pt, outer sep=0pt]
	\begin{scope}[shift={(100.0,100.0)}]
	\begin{scope}[draw=c888]
		\clip(-100,-100) rectangle (100,0);
		\path[draw]
				(-50.0000,-86.6025) -- (50.0000,-86.6025) (60.0000,-69.2820) -- (-60.0000,-69.2820)
				(-70.0000,-51.9615) -- (70.0000,-51.9615) (80.0000,-34.6410) -- (-80.0000,-34.6410)
				(-90.0000,-17.3205) -- (90.0000,-17.3205) (100.0000,0.0000) -- (-100.0000,0.0000)
				(-90.0000,17.3205) -- (90.0000,17.3205) (80.0000,34.6410) -- (-80.0000,34.6410)
				(-70.0000,51.9615) -- (70.0000,51.9615) (60.0000,69.2820) -- (-60.0000,69.2820)
				(-50.0000,86.6026) -- (50.0000,86.6026);
		\begin{scope}[rotate=120.0]
			\path[draw]
				(-50.0000,-86.6025) -- (50.0000,-86.6025) (60.0000,-69.2820) -- (-60.0000,-69.2820)
				(-70.0000,-51.9615) -- (70.0000,-51.9615) (80.0000,-34.6410) -- (-80.0000,-34.6410)
				(-90.0000,-17.3205) -- (90.0000,-17.3205) (100.0000,0.0000) -- (-100.0000,0.0000)
				(-90.0000,17.3205) -- (90.0000,17.3205) (80.0000,34.6410) -- (-80.0000,34.6410)
				(-70.0000,51.9615) -- (70.0000,51.9615) (60.0000,69.2820) -- (-60.0000,69.2820)
				(-50.0000,86.6026) -- (50.0000,86.6026);
		\end{scope}
		\begin{scope}[rotate=240.0]
			\path[draw]
				(-50.0000,-86.6025) -- (50.0000,-86.6025) (60.0000,-69.2820) -- (-60.0000,-69.2820)
				(-70.0000,-51.9615) -- (70.0000,-51.9615) (80.0000,-34.6410) -- (-80.0000,-34.6410)
				(-90.0000,-17.3205) -- (90.0000,-17.3205) (100.0000,0.0000) -- (-100.0000,0.0000)
				(-90.0000,17.3205) -- (90.0000,17.3205) (80.0000,34.6410) -- (-80.0000,34.6410)
				(-70.0000,51.9615) -- (70.0000,51.9615) (60.0000,69.2820) -- (-60.0000,69.2820)
				(-50.0000,86.6026) -- (50.0000,86.6026);
		\end{scope}
	\end{scope}
	\begin{scope}[draw=cccc]
		\clip(-100,0) rectangle (100,100);
		\path[draw] (0.0000,0.0000) circle (0.2822cm);
		\path[draw] (0.0000,0.0000) circle (0.8467cm);
		\path[draw] (0.0000,0.0000) circle (1.4111cm);
		\path[draw] (0.0000,0.0000) circle (1.9756cm);
		\path[draw] (0.0000,0.0000) circle (2.5400cm);
	\end{scope}
	\end{scope}

	\end{tikzpicture}

	\caption{for $\langle r^3\rangle$ symmetry, exactly $n$ small triangles must be in pairs which cross the diameter.}
\end{figure}

$\#c\langle r^2\rangle$ is A006366, ``number of cyclically symmetric plane partitions in the n-cube'' (Stanley's case 3), with formula $\#c\langle r^2\rangle(n) = \prod_{i=0}^{n-1} \frac{(3i+2)(3i)!(i+2n)!}{(2i+n)!(2i+n+1)!}$ \cite{Andrews79}.

$\#c\langle r\rangle$ is A049503, which is cryptically described as ``A005130$^2$'' and ``Expansion of generating function $A_{QT}^{(1)}(4n)$''. It corresponds to the cyclically symmetric and self-complementary plane partitions (Stanley's case 9). Note that again we must inflate the sequence (since $\langle r^3\rangle \le \langle r\rangle$ the same argument applies), and when we do so we get that $\#c\langle r\rangle(2m) = \prod_{i=0}^{m-1} \frac{(3i+1)!^2}{(i+m)!^2}$ \cite{Kuperberg94}.

We get the baseline
$$ \frac{1}{12} \prod_{i=0}^{n-1} \frac{i!(i+2n)!}{(i+n)!^2} + \frac{1}{4} \prod_{i=0}^{n-1} \frac{(2i)!(i+2n)!}{(2i+n)!(i+n)!} + \frac{1}{6} \prod_{i=0}^{n-1} \frac{(3i+2)(3i)!(i+2n)!}{(2i+n)!(2i+n+1)!} $$
and, for even $n = 2m$, the inflated correction
$$ \frac{1}{4} \binom{3m-1}{m} \prod_{i=0}^{2m-3} \frac{(i+4m)!(2i+2)!}{(2i+2m+2)!(i+2m)!} + \frac{1}{12} \prod_{i=0}^{m-1} \frac{i!^2(i+2m)!^2}{(i+m)!^4} + \frac{1}{6} \prod_{i=0}^{m-1} \frac{(3i+1)!^2}{(i+m)!^2} $$

\end{document}